\documentclass[11pt]{article}

\usepackage[latin1]{inputenc}
\usepackage{amsfonts,amsmath,amstext,amsbsy,amsthm,amscd,amssymb}
\usepackage{multirow}
\usepackage{graphicx}
\usepackage{epsf,epsfig}
\def\N{\mathbb{N}}
\def\R{\mathbb{R}}

\def\P{\mathbb{P}}
\def\E{\mathbb{E}}

\def\interior#1{\smash{\mathop{#1}\limits^{\lower1pt\hbox{$\scriptscriptstyle\circ$}}}}

\newtheorem{theo}{Theorem}[section]

\newtheorem{prop}{Proposition}[section]
\newtheorem{lemma}{Lemma}[section]

\title{{\bf  Computable infinite dimensional filters with applications
    to  discretized diffusion processes.}} 

\author{MIREILLE CHALEYAT-MAUREL$^{1}$ and VALENTINE GENON-CATALOT$^{2}$}

\date{\empty}

\begin{document}
\begin{titlepage}

\maketitle

\thispagestyle{empty}

{\sc \it \noindent
$^{1}$Universit\'e Ren\'e Descartes Paris 5, U.F.R. de Math\'ematiques
et Informatique, Laboratoire  MAP5 (CNRS UMR 8145) et
Laboratoire de Probabilit\'es et Mod\`eles Al\'eatoires (CNRS-UMR 7599), 45,
rue des Saints-P\`eres, 75270 Paris Cedex 06, France. e-mail: mcm@math-info.univ-paris5.fr\\
$^{2}$Universit\'e Ren\'e Descartes Paris 5, U.F.R. de Math\'ematiques
et Informatique, Laboratoire  MAP5 (CNRS-UMR 8145), 45, rue des Saints-P\`eres,
75270 Paris Cedex 06, France. e-mail: genon@math-info.univ-paris5.fr}

\begin{abstract}
Let us consider a pair signal-observation $((x_n,y_n), n \ge 0)$ where the
unobserved signal $(x_n)$ is a Markov chain and the observed component
is such that, given the whole sequence $(x_n)$, the random variables
$(y_n)$ are independent and the conditional distribution of $y_n$ only
depends on the corresponding state variable $x_n$. The main problems raised by
these observations are the prediction and filtering of $(x_n)$. We
introduce sufficient conditions allowing to obtain computable filters using mixtures of
distributions. The filter system may be finite or infinite dimensional. The method is applied to the case where the signal
$x_n=X_{n \Delta}$ is a discrete sampling of a one dimensional
diffusion process: Concrete models are proved to fit in our
conditions. Moreover, for these models, exact likelihood inference
based on the observation $(y_0, \ldots, y_n)$ is
feasable.   
\end{abstract}

\vspace{1cm}

\vspace{1cm}

\vspace{0,5cm}
\noindent
{\bf Running title:} Computable filters.
\end{titlepage}

\section{Introduction}

Let us consider a pair signal-observation $((x_n,y_n), n \ge 0)$ where the
unobserved signal $(x_n)$ is a Markov chain and the observed component
is such that, given the whole sequence $(x_n)$, the random variables
$(y_n)$ are independent and the conditional distribution of $y_n$ only
depends on the corresponding state variable $x_n$. This is a
classical setting in the field of  non linear filtering and  the process $(y_n)$
is often called a hidden Markov model.

 In this context, a  central problem  that has been the subject of a
 huge number of contributions is the study of the exact filter, {\em
   i.e.} the sequence of conditional distributions of $x_n$ given
 $y_n, \ldots, y_{1}, y_{0}$, $n \ge 0$. On the other hand,
 statistical inference based on the (non Markovian) observations
 $(y_0, \ldots, y_n)$ requires the knowledge of the successive
 conditional distributions   of $y_n$ given $y_{n-1}, \ldots,
 y_{0}$. These are obtained through the prediction filter, {\em i.e.}
 the sequence of conditional distributions of $x_n$ given $y_{n-1}, \ldots, y_{1}, y_{0}$.
Although the exact and the prediction filter may both be calculated
recursively by an explicit algorithm,  iterations become rapidly
intractable and exact formulae are difficult to
obtain. To overcome this difficulty, authors generally try to find a
parametric family ${\cal F}$ of distributions on the state space
${\cal X}$ of $(x_n)$ ({\em i.e.} a family of distributions specified
by a finite fixed number of real parameters) such that if ${\cal L}(x_0) \in {\cal F}$, then,
for all $n$, ${\cal L}(x_n|y_n, \ldots, y_{1}, y_{0})$ and ${\cal
  L}(x_n|y_{n-1}, \ldots, y_{1}, y_{0})$ both belong to ${\cal F}$. In
this case, the model $((x_n,y_n), {\cal F})$ is called a
finite-dimensional filter system and it is enough to describe each
conditional distribution by the parameters that characterize it. This situation
is illustrated by the linear Gaussian Kalman filter (see below Section \ref{kalm}). Whenever the initial
distribution of the signal is Gaussian, specified by its mean and variance, then,
all the successive conditional distributions are Gaussian and there is
an explicit recursive algorithm which gives the stochastic process of
the conditional  means and
variances.

Necessary and sufficient conditions for the existence of
finite-dimensional filters in discrete time have been given in Sawitzki (1981)
(see also Runggaldier and Spizzichino (2001)). The case of continuous
time filters was treated in Chaleyat-Maurel and Michel (1984). As a
consequence of these papers, it appears that very few finite
dimensional filters are available and they are often obtained as the
result of an {\em ad hoc} construction (see {\em e.g.}, the new
constructive approach presented in Ferrante and Vidoni (1998) and the
references therein).

In what follows, we propose a method to obtain computable
filters that may be finite or infinite dimensional. Our approach is a
generalization of the one developped in Genon-Catalot (2003) and
Genon-Catalot and Kessler (2004) for a special model. The method
is well fitted   for the filtering of discretized diffusion
processes and illustrated with examples. 

More precisely, (Section \ref{suff})
we  consider  at first a sequence of parametric families of distributions
${\cal F}^{i}= \{\nu^{i}_{\theta}, \theta \in \Theta \}$, $i \in
\N$, where $\Theta \subset \R^{p}$ is a parameter set. Then, we
construct an enlarged family  by means of mixtures. Let us define the
set $S$ of mixture parameters:
 \begin{equation} \label{mixt}
S= \{\alpha=(\alpha_i, i \ge 0), \forall i \ge 0, \alpha_i \ge 0, \sum_{i=0}^{\infty} \alpha_i=1\}.
\end{equation}
Then, we set:
\begin{equation}
{\bar {\cal F}}= \{\nu=\sum_{i \ge 0} \alpha_i \nu^{i}_{\theta}=\nu_{\theta,\alpha},\;
\alpha=(\alpha_i, i \ge 0) \in S,\;  \theta \in \Theta \}.
\end{equation} 
Each distribution $\nu=\nu_{\theta,\alpha}$  in the above class is
thus specified by an usual
parameter $\theta$ and a mixture parameter $\alpha=(\alpha_i, i \ge
0)$. We give sufficient conditions on the class ${\cal F}=
\cup_{i \in \N} {\cal F}^{i}$  ensuring that if $\nu={\cal L}(x_0)$
belongs to ${\bar {\cal F}}$, then, the exact and the prediction
filter both evolve within ${\bar {\cal F}}$. These conditions involve
the conditional distribution of $y_n$ given $x_n=x$ and the transition
operator of the hidden chain $(x_n)$. Of course, the most interesting case
is when the mixture distributions obtained for the filters have a
finite number of components. To this end, we introduce the sub-class
${\bar {\cal F}}_{f}$ of ${\bar {\cal F}}$ composed of distributions
$\nu=\nu_{\theta,\alpha}$ such that $\alpha_i=0$ for $i$ greater than
some integer $N$. We give sufficient conditions ensuring that, when
${\cal L}(x_0) \in {\bar {\cal F}}_{f}$, then, the exact and the
prediction filter evolve in ${\bar {\cal F}}_{f}$ (Conditions
(C1)-(C2) and Theorem \ref{cs}). They are thus specified by a finite
number of parameters. This finite number may change along the
iterations. Nevertheless, the filters are explicitly and exactly
computable. Let us note that our method has links with the one
developped in Di Masi et al (1983). In this paper, the distributions
of the filters are allowed to be finite linear combinations of
parametric distributions. The number of terms in the linear
combination may also  change along the iterations. However, the
coefficients in each linear combination are possibly negative. So
these distributions are not mixtures of parametric distributions and
their interpretation is therefore difficult.

In Section \ref{kalm}, we illustrate our method on the classical
Kalman filter with non Gaussian initial distributions. We introduce an appropriate class of non Gaussian
distributions and give the corresponding  explicit formulae for the
filters. In Sections \ref{scalepert}-\ref{scalecir}-\ref{poisson}, we
consider models satisfying our sufficient conditions and for which the
signal $x_n=X_{n \Delta}$ is a discrete sampling of a one-dimensional
diffusion process. In Section \ref{scalepert}, we study the
observation equation $y_n=x_n w_n$  where the signal $(x_n)$ and the
noise $(w_n)$ are independent sequences of positive random variables. The signal $x_n$ is a discrete
sampling of the diffusion process given by the stochastic differential
equation
\begin{equation} 
dX_t= (\theta X_t+ \frac{(\delta-1) \sigma^{2}}{2 X_t}) dt +  \sigma  d\beta_t,
\end{equation}
where $(\beta_t)$ is a standard one-dimensional Brownian motion,
$\theta$ is a real parameter, $\sigma$ is positive and $\delta$ is a
real number such that $\delta >1$. This process is called the radial
Ornstein-Uhlenbeck process. In the case $\delta=1$, we consider the
absolute value of a one dimensional Ornstein-Uhlenbeck
process. The noises $w_n$ have a specific distribution to ensure
explicit formulae. In Sections \ref{scalecir} and \ref{poisson}, we
exploit the results of Section \ref{scalepert} to  develop models
based on the same signal but different observation equations. In particular, in Section \ref{scalecir}, we propose some stochastic
volalility type models. Section
\ref{concl} contains some concluding remarks. Some technical proofs
are given in the Appendix.


\section{Sufficient conditions for computable 
  filters.} \label{suff}
\subsection{The filtering- prediction algorithm.}

We introduce our notations in an abstract framework that will become
concrete through the examples below.

We denote by ${\cal X}$ the state-space of the unobserved Markov chain
$(x_n)$, which is equipped with a sigma-field ${\cal B}$. We assume
that this chain is time-homogeneous and denote  by $P$ its
transition operator. The time-homogeneity assumption is  here for the
sake of simplicity of notations. It is not essential
 and may be relaxed (see Section \ref{concl}).

 The state-space of the observed component  $(y_n)$ will be
denoted by ${\cal Y}$ and its sigma-field by ${\cal C}$. We assume that the conditional distribution  of
$y_n$ given $x_n=x$ does not depend on $n$ and that this distribution
is given by a density with respect to a common dominating measure
$\mu$ on ${\cal Y}$:
\begin{equation} \label{condit}
F_{x}(dy)={\cal L}(y_n|x_n=x)= f_{x}(y) \mu(dy).
\end{equation}
As stated above, this is also the conditional distribution of $y_n$
given the whole sequence $(x_k)$ when $x_n=x$. Again, the
time-homogeneity assumption here is not essential. On the contrary, the
existence of a common dominating measure $\mu$ for the family of distributions
$(F_{x}(dy), x \in {\cal X})$ is essential for the filtering-prediction
recursive equations. In the statistical vocabulary, it means that this
family of distributions is a dominated family when $x$ is considered
as a parameter.

In the concrete models that we investigate in further sections, we shall take ${\cal
  X}$  equal to $\R$ or $(0, + \infty)$, and  ${\cal Y}$ equal to
$\R$, $(0,+\infty)$ or $\N$.

Let us now briefly recall the filtering-prediction algorithm (see {\em e.g.} Del
Moral and Guionnet (2001)). Given an initial distribution for $x_0$, there is a well known algorithm that allows to compute the successive conditional distributions:
\begin{eqnarray*}
{\cal L}(x_0)\xrightarrow[]{{\rm  updating}}  {\cal L}(x_0\vert y_0)
\xrightarrow[]{{\rm prediction}} {\cal L}(x_1\vert y_0)\hphantom{uuuuuuuuuuuuuu} \\ \hphantom{uuuuuuuuuuuuuu}\xrightarrow[]{{\rm updating}} {\cal L}(x_1\vert y_1,y_0)\xrightarrow[]{{\rm prediction}} {\cal L}(x_2\vert y_1,y_0)\ldots.
\end{eqnarray*}

By the above iterations, we get two kinds of  distributions on ${\cal X}$:
\begin{equation} \label{opt}
\nu_{n|n:0}= {\cal L}(x_n|y_n, \ldots, y_0),
\end{equation}
\begin{equation} \label{predic}
\nu_{n+1|n:0}= {\cal L}(x_{n+1}|y_n, \ldots, y_0),
\end{equation}
The distribution (\ref{opt}) is called the optimal or exact filter and
(\ref{predic}) is the prediction filter.

These are obtained using two steps: the updating and the prediction
steps which can be described by introducing the following
operators. Let ${\cal P}({\cal X})$ denote the set of probability measures
on ${\cal X}$. For $\nu \in {\cal P}({\cal X})$, the probability
$\varphi_{y}(\nu)$ is defined by (see (\ref{condit})):
\begin{equation} \label{phi}
\varphi_{y}(\nu)(dx)= \frac{f_{x}(y) \nu(dx)}{p_{\nu}(y)},
\end{equation}
(with the convention that $0/0=0$), where the denominator  is equal to
\begin{equation} \label{marg}
p_{\nu}(y)=  \int_{\cal X}\nu(d\xi) f_{\xi}(y).
\end{equation}
The operator $\varphi_y$ is called the up-dating operator which
allows to take into account a new observation. 

On the other hand,  the prediction operator is as follows. For $\nu \in
{\cal P}({\cal X})$, the probability measure 
\begin{equation} \label{psiP}
\psi(\nu)=\nu P
\end{equation}
is
obtained by applying the transition operator of the hidden Markov
chain. It is 
defined by:
\begin{equation} \label{psiP2}
\psi(\nu)(dx')= \nu P(dx')= \int_{{\cal X} }
 \nu(dx) P(x,dx').
\end{equation}

We have, for $n \ge 0$, (with $\nu_{0|-1:0}=\nu_0={\cal L}(x_0))$
\begin{equation}
\nu_{n|n:0}=\varphi_{y_n}(\nu_{n|n-1:0}), \quad \nu_{n+1|n:0}=\psi(\nu_{n|n:0}).
\end{equation}
The optimal filter  is obtained via the operator
\begin{equation} \label{Phihat}
{\hat \Phi}_{y}=\varphi_{y} \circ \psi,
\end{equation}
and 
\begin{equation}
\nu_{n+1|n+1:0}={\hat \Phi}_{y_{n+1}}(\nu_{n|n:0})= {\hat
  \Phi}_{y_{n+1}} \circ \ldots \circ {\hat \Phi}_{y_{1}} \circ \varphi_{y_0}(\nu_0).
\end{equation}

 The prediction filter is obtained via the operator
\begin{equation} \label{Phi}
\Phi_{y}=\psi \circ \varphi_{y} ,
\end{equation}
and 
\begin{equation}
\nu_{n+1|n:0}=\Phi_{y_{n}}(\nu_{n|n-1:0})=\Phi_{y_{n}} \circ \ldots
\circ \Phi_{y_{0}}(\nu_0) .
\end{equation}
Moreover, the conditional density of $y_n$ given $(y_{n-1}, \ldots,
y_0)$ is obtained as the following marginal density (see (\ref{marg})):
\begin{equation} \label{margn}
p(y|y_{n-1}, \ldots, y_0)= p_{\nu_{n|n-1:0}}(y).
\end{equation}
And the exact density of $(y_0,y_1, \ldots, y_n)$ is obtained as the
product of the successive conditional densities $p(y_i|y_{i-1}, \ldots, y_0)$.
\subsection{Sufficient conditions.}
The iterations above are rapidly untractable unless both
operators  (\ref{phi}) and (\ref{psiP}) evolve in a parametric family of
distributions, {\em i.e.}  distributions specified by a fixed
finite number  of real parameters. In what follows, this number of
parameters will possibly vary along iterations.

More precisely, let us define a class
${\bar {\cal F}}$ of distributions on ${\cal X}$ as follows. First, we start
with a parametric  class of the form
\begin{equation} \label{f}
{\cal F}=  \{\nu_{\theta}^{i}, \theta \in
\Theta, i \in \N \}
\end{equation}
  where $\Theta \subset \R^{p}$ is a parameter set. Each
distribution in  ${\cal F}$ is thus specified by a couple
$(i,\theta) \in \N \times \Theta$ and there is a one-to-one
correspondence between  $\N \times \Theta$
and the class ${\cal F}$. Then, using  the set S of mixture
parameters  defined in (\ref{mixt}), we build the enlarged
class composed of  convex combinations of distributions
$\nu^{i}_{\theta}$ having the same parameter $\theta$:  
\begin{equation} \label{fbar}
{\bar {\cal F}}= \{\nu; \nu= \sum_{i=0}^{\infty}
\alpha_i\;\nu_{\theta}^{i}=\nu_{\theta, \alpha}, \theta \in \Theta,
\alpha= (\alpha_i, i \ge 0) \in S\}.
\end{equation}
Now, each distribution $\nu_{\theta, \alpha}$ on ${\bar {\cal F}}$ is
specified by a  parameter $\theta$ and a mixture parameter
$\alpha$. We stress the fact that all components in a given mixture have the same parameter
$\theta$. The mixture parameter $\alpha$  of a $\nu_{\theta, \alpha}$
may  or may not depend on $\theta$. For
\begin{equation}
\alpha= \alpha^{(i)} \quad \mbox{given by} \quad
\alpha_i^{(i)}=1, \quad \alpha_j^{(i)}=0, j \neq i,
\end{equation}
we get the distribution $\nu_{\theta}^{i}$:
 \begin{equation} \label{nui}
\nu_{\theta}^{i}=\nu_{\theta,\alpha^{(i)}}.
\end{equation}
Obviously, ${\cal F} \subset {\bar {\cal F}}$. But the resulting
extended class may be considerably larger. Of course,  the number of
components in the mixture can be finite. So, we shall define the
length of  a
mixture parameter $\alpha$ by
\begin{equation} \label{l}
l(\alpha)= \sup{\{i; \alpha_{i} >0\}}.
\end{equation}
We define the sub-class of distributions with finite-length mixture parameter by
\begin{equation} \label{ff}
{\bar {\cal F}}_{f}= \{ \nu= \nu_{\theta, \alpha} \in {\bar {\cal F}},
l(\alpha)<\infty \}.
\end{equation}
On the other hand, when $\alpha$ has infinite length, the series
defining an element $\nu$ in (\ref{fbar}) may have an explicit sum, which
will be another expression of  $\nu$.  

Now, we want conditions such that, for $\nu \in {\bar {\cal
    F}}$,  $\varphi_y(\nu)$ and
$\psi(\nu)$ both belong to ${\bar {\cal F}}$. In such a case, it will
be enough to express both operators in terms of the couple $(\theta,\alpha)$ specifying the distributions in  ${\bar {\cal
    F}}$. Moreover, when the two operators evolve within ${\bar {\cal
    F}}_{f}$, then the exact and the prediction filters are exactly
computable even if the number of mixture components varies along the iterations.

Let us  consider the following conditions.
\begin{itemize}
\item [$\bullet$ (C1)] For all $y \in {\cal Y}$,  for all $\nu \in
  {\cal F}$, $\varphi_y(\nu) \in {\cal F}$ (see (\ref{phi})). More
  precisely, for all $(i,\theta) \in \N \times \Theta$, 
$$
\varphi_y(\nu_{\theta}^{i}) = \nu_{T_{y}(\theta)}^{t_{y}(i)},
$$
with $T_{y}(\theta) \in \Theta$, $t_{y}: \N
\rightarrow \N$  a one-to-one mapping  and $(\theta,y)
\rightarrow T_y(\theta)$  measurable. 
\item [$\bullet$ (C2)] For all $\nu \in {\cal F}$, $\psi(\nu)=\nu P
  \in {\bar {\cal F}}$. More precisely, for all $(i,\theta)$,
  $\psi(\nu_{\theta}^{i})$  may be written as
$$
\psi(\nu_{\theta}^{i})= \sum_{j \ge 0} \alpha_{j}^{(i,\theta)} \nu_{\tau(\theta)}^{j},
$$
where $\alpha^{(i,\theta)} \in S$,  $\tau(\theta)
\in \Theta$ and $\theta \rightarrow (\tau(\theta),\alpha^{(i,\theta)})
$ measurable.
\item [$\bullet$ (C2-f)] For all $\nu \in {\cal F}$, $\psi(\nu)=\nu P
  \in {\bar {\cal F}_f}$, with, using the notations of (C2),  for all $(i,\theta)$, $l(\alpha^{(i,\theta)})=L(i) < \infty$ and
  the mapping $i \rightarrow L(i)$ is non decreasing.
\item [$\bullet$ (C3)] For all $x \in {\cal X}$, $P(x,dx')$ belongs to the class
  ${\bar {\cal F}}$, and may be written as
$$
P(x,dx')= \sum_{i \ge 0} \alpha_{i}^{0}(x) \nu_{\theta_0}^{i}(dx'),
$$ 
where $\alpha^{0}(x) \in S$, $\theta_0 \in \Theta$ and $x \rightarrow
\alpha^{0}(x)$ is measurable.
\end{itemize}

Let us make some comments about these conditions. Condition (C1)
concerns only $\varphi_y$ and the class ${\cal F}$. The up-dating
operator $\varphi_{y}$ becomes the following mapping from
${\cal F}$ to ${\cal F}$:
\begin{equation} \label{nufib}
(i,\theta) \rightarrow (t_{y}(i), T_{y}(\theta)).
\end{equation}
Note that, in condition (C1), the function $t_{y}(.)$ must not depend on $\theta$.
 Thus,  the class  ${\cal F}$ is a conjugate class of distributions
for the parametric family $F_{x}(dy)= f_{x}(y) d\mu(y)$ (in the sense
of Bayesian estimation). Conditions (C2)-(C2-f)-(C3)  concern the
transition operator $P$ and the class ${\cal F}$. Condition (C3) 
implies that, when the signal starts at a fixed  $x_0=x$, then,  the
distribution of $x_1$  belongs to the enlarged class ${\bar {\cal
  F}}$. Therefore, we can consider that Dirac measures belong to the
enlarged class or directly add all Dirac measures to this
class. Conditions (C2-f) and (C3) may appear contradictory. Actually,
this is not the case because a distribution $\nu$ in ${\bar {\cal F}}$ may
have two different representations, {\em i.e.} the equality  $\nu_{\theta, \alpha}=
\nu_{\theta', \alpha'}$ does not imply $(\theta, \alpha)=(\theta',
\alpha')$. Moreover, one representation may be finite and the other
infinite. We discuss this point in Section \ref{scalepert}. 

We have the following result.
\begin{theo} \label{cs}
\begin{enumerate}
\item Assume (C1)-(C2). If $\nu$ belongs to ${\bar {\cal F}}$
  (see (\ref{fbar})),
  then, $\varphi_{y}(\nu)$ and $\psi(\nu)$ both belong to ${\bar {\cal
      F}}$.
 More precisely, if $\nu=\nu_{\theta,\alpha}$ then,
  $$\varphi_{y}(\nu)=\nu_{T_{y}(\theta),a_y(\theta,\alpha)} \quad
  \mbox{and} \quad
  \psi(\nu)=\nu_{\tau(\theta),b(\theta,\alpha)}$$ 
where
  $T_y(\theta), \tau(\theta)$ are defined in Conditions (C1)-(C2) and
  the mixture parameters 
  $a_{y}(\theta, \alpha)$ and $b(\theta, \alpha)$  are given in
  formulae (\ref{nufif}) and (\ref{psialf}) or (\ref{psif}). 
\item Assume (C1)-(C2-f). If $\nu$ belongs to ${\bar {\cal F}}_{f}$ (see
  (\ref{ff})),
  then, $\varphi_{y}(\nu)$ and $\psi(\nu)$ both belong to ${\bar {\cal
    F}}_{f}$. 
\end{enumerate}
\end{theo}
\proof Consider first the up-dating operator $\varphi_{y}$ (see (\ref{condit})-(\ref{phi})-(\ref{marg})).  Let
$\nu=\nu_{\theta,\alpha}= \sum_{i \ge 0} \alpha_i \;
\nu_{\theta}^{i}$. Then (see (\ref{marg})), 
\begin{equation} \label{margi}
p_{\nu}(y)= \sum_{i \ge 0} \alpha_i \;
p_{\nu_{\theta}^{i}}(y),  
\end{equation}
with, for all $i$, 
\begin{equation}
p_{\nu_{\theta}^{i}}(y)= \int_{{\cal X}} \nu_{\theta}^{i}(dx) f_x(y). 
\end{equation}
We have
\begin{equation}
f_x(y)\nu(dx)= \sum_{i \ge 0} \alpha_i \;\nu_{\theta}^{i}(dx) f_x(y).
\end{equation} 
Using (C1), since the mapping $t_{y}$  is one-to-one, we get
\begin{equation}
f_x(y)\nu(dx)= \sum_{i \ge 0} \alpha_i \;p_{\nu_{\theta}^{i}}(y)
\nu_{T_y(\theta)}^{t_{y}(i)}(dx)=
\sum_{j \ge 0} \alpha_{t_{y}^{-1}(j)} \;p_{\nu_{\theta}^{t^{-1}(j)}}(y) \nu_{T_y(\theta)}^{j}(dx)
\end{equation}
Hence, 
\begin{equation}
\varphi_{y}(\nu_{\theta, \alpha)}= \nu_{T_{y}(\theta),a_{y}(\theta,\alpha)},
\end{equation}
where the parameter $T_{y}(\theta)$ is defined in (C1) and the mixture
coefficient $a_{y}(\theta,\alpha)$ is given by (see (\ref{margi})):
\begin{equation} \label{nufif}
a_{y}(\theta, \alpha)_{j}= \frac{\alpha_{t_{y}^{-1}(j)}\,p_{\nu_{\theta}^{t_{y}^{-1}(j)}}(y)}{p_{\nu}(y)}.
\end{equation}
The operator $\varphi_{y}$ on ${\bar {\cal F}}$ can be expressed in
terms of the parameters  by the
mapping:
\begin{equation}
(\theta, \alpha) \in \Theta \times S \rightarrow (T_{y}(\theta),
a_{y}(\theta,\alpha)) \in \Theta \times S.
\end{equation}

By (\ref{nufif}), the mixture parameter $a_{y}(\theta,\alpha)$ has
finite length when $\alpha$ has finite length.

Consider now the prediction operator $\psi$ (see (\ref{psiP})). By linearity and (C2),
we get
$$
\psi(\nu)= \sum_{i \ge 0} \alpha_i \; \psi(\nu_{\theta}^{i})=\sum_{j \ge 0} b(\theta,\alpha)_j \; \nu_{\tau(\theta)}^{j}=
\nu_{\tau(\theta), b(\theta,\alpha)},
$$
where $\tau(\theta)$ defined in (C2) and the new mixture parameter is
obtained by interchanging sums and is 
given by
\begin{equation} \label{psialf}
b(\theta,\alpha)_j= \sum_{i \ge 0} \alpha_i \; \alpha_{j}^{(i,\theta)}.
\end{equation}
The prediction operator $\psi$ is therefore now given defined  by
the mapping:
\begin{equation}
(\theta, \alpha) \rightarrow (\tau(\theta), b(\theta,\alpha)).
\end{equation}
Now, if $l(\alpha)=p$, and (C2-f) holds, then
\begin{equation} \label{psif}
b(\theta,\alpha)_j= \sum_{i=0}^{p} \alpha_i\;\alpha_j^{(i,\theta)}\;1_{(j
  \le L(i))} = \sum_{i \ge L^{-1}(j),i \le p} \alpha_i \alpha_{j}^{(i,\theta)}
\end{equation}
where $L^{-1}(j)= \inf\{i;L(i) \ge j\}$.
Since $\alpha_i=0$ for $i >p$, $b(\theta,\alpha)_j=0$ as soon as
$L^{-1}(j)>p$.
\endproof

\noindent
{\bf Remark.}
\begin{enumerate}
\item It is worth noting that our conditions imply that the
the parameters $T_y(\theta), \tau(\theta)$ only depend on $\theta$
whereas $a_y(\theta, \alpha), b(\theta, \alpha)$  depend on $(\theta,
\alpha)$.
\item For  $ x \in
  {\cal X}$, it is immediate to check that
  $\varphi_y(\delta_x)=\delta_x$. By (C3), $\psi(\delta_x)=
  \nu_{\theta_{0}, \alpha^{0}(x)}$. So the algorithm starting with a
  deterministic initial condition evolves in ${\bar {\cal F}}$.  
\end{enumerate}
\section{The Kalman filter with non Gaussian initial condition.}
\label{kalm}
 Our first example is based on  the classical and simplest standard one-dimensional
Kalman filter. It is well
known (see {\em e.g.} Makowski (1986)) that, whatever the initial
distribution for the Kalman filter, it is possible to compute the
prediction and exact filters. We illustrate this property through a
special family of initial distributions. Let us recall the model. The
observation equation is   given by
\begin{equation} \label{kalman}
y_n= h\, x_n + \gamma \,w_n,
\end{equation}
with $h, \gamma$ constants ($\gamma>0$), $(w_n)$ a standard one-dimensional 
Gaussian white noise. And for the signal
\begin{equation} \label{ar}
x_n= a\, x_{n-1} + \beta \,\eta_n,
\end{equation}
with $a, \beta$ constants ($\beta>0$), $(\eta_n)$ a standard
one-dimensional Gaussian white noise. The sequences $(x_n)$ and
$(w_n)$ are assumed to be independent. Now, the conditional distribution of $y_n$ given
$x_n=x$ is
\begin{equation} \label{loicond}
f_x(y)dy= {\cal N}(h\,x, \gamma^{2}).
\end{equation}
And the transition operator of $(x_n)$ is
 \begin{equation} \label{transop}
P(x,dx')=p(x,x')dx'= {\cal N}(ax, \beta^{2}).
\end{equation}
We introduce below a class of non Gaussian distributions and show
that our conditions (C1)-(C3) hold for this class. Therefore,
(\ref{opt})-(\ref{predic}) can be explicitely computed. Before doing
this, we recall the classical case.

  \subsection{The standard Kalman filter.}
It is well-known that if the initial distribution is Gaussian (or deterministic) then,
for all $n$, the distributions (\ref{opt}) and (\ref{predic}) are
Gaussian. 
Let us denote by ${\cal G}= \{ {\cal N}(m, \sigma^{2}), m \in \R,
\sigma^{2}> 0\}$ the class
of Gaussian distributions. The up-dating and prediction operators are
from ${\cal G}$ onto ${\cal G}$.  And, some classical and elementary computations yield:
\begin{itemize}
\item The up-dating step is: $\nu={\cal N}(m, \sigma^{2}\} \rightarrow \varphi_y(\nu)={\cal
    N}({\hat m}, {\hat \sigma}^{2})$, with 
\begin{equation} \label{chap}
{\hat m}= \frac{m \gamma^{2}+hy\sigma^{2}}{\gamma^{2}+h^{2}\sigma^{2}},\quad  {\hat \sigma}^{2}=\frac{\sigma^{2}\gamma^{2}}{\gamma^{2}+h^{2}\sigma^{2}}.
\end{equation}
\item The marginal distribution is: $p_{\nu}(y)dy = {\cal N}(hm, \gamma^{2}+h^{2}\sigma^{2})$
\item The prediction step is: 
$\nu={\cal N}( m,  \sigma^{2}) \rightarrow
  \psi(\nu)={\cal N}({\bar m}, {\bar \sigma}^{2})$ with
\begin{equation} \label{bar}
{\bar m}= a\, m, \quad {\bar \sigma}^{2}=\beta^{2}+ a^{2} 
    \sigma^{2}.
\end{equation}

\end{itemize}
The formulae above hold true when $\sigma=0$ 
allowing to include the case of Dirac measures.
Note that, for all $x$,  $P(x,dx')= {\cal N}(ax, \beta^{2})$ also
belongs to ${\cal G}$. 
\subsection{An extended Kalman filter.}
Now, we enlarge the class of Gaussian distributions using new
distributions and  mixtures. Consider three parameters $\mu, m, \sigma^{2}$ with $\mu,m \in \R$ and
$\sigma^{2}>0$. For $i=0$, set
\begin{equation} \label{gnu0}
\nu^{0}_{(0,m,\sigma^{2})}(dx)= {\cal N}(m,\sigma^{2}).
\end{equation}
and ${\cal F}^0= {\cal G}$.
For $i \ge 1$, set
\begin{equation} \label{gnui}
\nu^{i}_{(\mu,m,\sigma^{2})}(dx)= \frac{(x+\mu)^{2i}}{C_{2i}(m+\mu;\sigma^{2})}\nu^{0}_{(0,m,\sigma^{2})}(dx),
\end{equation}
where the normalizing constant is given by
\begin{equation} \label{cst1}
C_{2i}(m+\mu;\sigma^{2})= \E((\sigma X+ \mu+m)^{2i}),
\end{equation}
for $X$ a standard Gaussian random variable. Let  denote by
\begin{equation} \label{momg}
C_{2i}= \E(X^{2i})= \frac{(2i)!}{2^{i}i!},
\end{equation}
the $2i$-th moment of $X$. Then, for $0 \le k \le i$,  some elementary computations yield:
\begin{equation} \label{formule}
\binom{2i}{2k} C_{2(i-k)}= \frac{C_{2i}}{C_{2k}} \binom{i}{k}.
\end{equation}
We deduce
\begin{equation} \label{cst2}
C_{2i}(m+\mu;\sigma^{2})= \sum_{k=0}^{i} (m+\mu)^{2k} \sigma^{2(i-k)}
\frac{C_{2i}}{C_{2k}} \binom{i}{k}.
\end{equation}
Set ${\cal F}^i= \{\nu^{i}_{(\mu,m,\sigma^{2})},  (\mu,m,\sigma^{2})
\in \R \times \R \times (0,+\infty)\}$.
Define ${\cal F}= \bigcup_{i \ge 0} {\cal F}^i$. 

In the Appendix, we study some elementary properties of  these
distributions.  

The class ${\bar {\cal F}}$ is defined as in (\ref{fbar}). All distributions in
${\bar {\cal F}}$ have density with respect to a Gaussian law. 

Now, we check conditions (C1)-(C2-f)-(C3). Condition (C3) evidently
holds (see (\ref{transop})) since ${\bar {\cal F}}$ contains all Gaussian
distributions. By the following proposition, condition (C2-f) holds.
\begin{prop} \label{upd}
 Consider the model given by
  (\ref{kalman})-(\ref{ar})-(\ref{loicond})-(\ref{transop}). 
\begin{enumerate}
\item For $i=0$, and $\nu= {\cal N}(m, \sigma^{2})$,  $\varphi_y(\nu)= {\cal
    N}({\hat m}, {\hat \sigma}^{2})$ with ${\hat m}, {\hat \sigma}^{2}$
    given in (\ref{chap}).  
\item For $i \ge 1$, and  $\nu=\nu_{(\mu,m,\sigma^{2})}^{i}(dx)$ in  ${\cal F}^i$, $\varphi_y(\nu)=\nu_{(\mu,{\hat m},{\hat
      \sigma}^{2})}^{i}$. 
\end{enumerate}
Condition (C1) holds with $t_{y}(i)=i$, for all $i \ge 0$ and
$T_y(\mu,m,\sigma^{2})= (\mu,{\hat m}, {\hat \sigma}^{2})$ (where for
$i=0$, $\mu=0$). 
\end{prop}
The proof is obtained in the same way as for the up-dating step for the
classical Kalman filter (see the Appendix).
Now, looking at formulae (\ref{marg}) and (\ref{psiP2}), since $f_x(y)$
and $p(x,x')$ in this model are both Gaussian kernels, the computation
of marginal distributions and the checking of (C2-f) are
identical up to a change of notations. The results are given in the
following proposition.
\begin{prop} \label{pre}
 Let $\nu=\nu_{(\mu,m,\sigma^{2})}^{i}(dx)$ belong to
  ${\cal F}^i$.
\begin{enumerate}
\item Then 
\begin{equation} \label{psinu}
\psi(\nu)= \sum_{k=0}^{i} {\bar \alpha}_k^{(i)} \nu_{({\bar
    \mu},{\bar m}, {\bar \sigma}^{2})}^{k},
\end{equation}
with ${\bar m}, {\bar \sigma}^{2}$ given in (\ref{bar}),
\begin{equation} \label{mubar}
{\bar \mu}=\frac{m \beta^{2}+ \mu {\bar \sigma}^{2}}{a \, \sigma^{2}},
\end{equation}
and for $k=0, \ldots, i$,
\begin{equation} \label{albar}
{\bar \alpha}_k^{(i)}=  \binom{i}{k}
\frac{\beta^{2(i-k)}}{B_{i}}\sum_{j=0}^{k}
\binom{k}{j}\frac{(\mu+m)^{2j}}{C_{2j} \sigma^{2j}}\frac{a^{2(k-j)}
  \sigma^{2(k-j)}}{{\bar \sigma}^{2(i-j)}}, 
\end{equation}
with 
\begin{equation} \label{bi}
B_i= \sum_{k=0}^{i} \binom{i}{k} \frac{(\mu+m)^{2k}}{C_{2k} \sigma^{2k}}
\end{equation}
($\sum_{k=0}^{i} {\bar \alpha}_k^{(i)}=1$).
Hence, Condition (C2-f) holds with $\tau(\mu,m,\sigma^{2})= ({\bar \mu},
{\bar m}, {\bar \sigma}^{2})$, for $k=0, \ldots,i$, $\alpha_k^{(i,\mu,m,\sigma^{2})}={\bar \alpha}_k^{(i)} $ and the length of the mixture
parameter of $\psi(\nu)$ for $\nu \in {\cal F}^i$
  is $L(i)=i$. 
\item The marginal distribution $p_{\nu}(y)dy$ is given by  the same
  formula as (\ref{psinu}) with $(a,\beta^{2})$ everywhere replaced by $(h,\gamma^{2})$.
\end{enumerate}
\end{prop}
The proof is given in the Appendix. Note that, using $\frac{a^2
  \sigma^{2}}{{\bar \sigma}^{2}}=
1-\frac{\beta^2}{{\bar \sigma}^{2}}$, we have 
$$
{\bar m}+{\bar \mu}= a(m+\mu)\left( 1-\frac{\beta^2}{{\bar \sigma}^{2}}\right)^{-1}.
$$
\noindent
{\bf Remark.} 
It is worth noting that, in this model, the number of parameters
remains fixed along iterations: If the initial condition is specified
by parameters $(\mu,m, \sigma, \alpha)$ with $l(\alpha)=p$, {\em i.e.}
$3+p+1$ parameters, then the length of the mixture parameter will
always be equal to $p$ and the number of parameters will remained
fixed equal to $3+p+1$. This is not surprising since the Kalman
filter is a finite-dimensional filter, even when the initial condition
is non Gaussian (see {\em e.g.} Makowski (1986)).
\section{Scale perturbation of a radial Ornstein-Uhlenbeck process.} \label{scalepert}

In this section, we consider multiplicative perturbation models of the form
\begin{equation} \label{mult}
y_n=  x_n w_n
\end{equation}
where $(w_n)$ is a sequence of i.i.d. positive random variables and
$(x_n)$ is also a positive signal independent of the sequence $(w_n)$. The
multiplicative structure comes from the field of Finance with the
so-called stochastic volatility models. However, in stochastic
volatility models, the noises are standard Gaussian variables (see the
next section). The advantage of positive signal and noise is that we
can interpret the model as a scale perturbation of a positive signal.

Now, we consider a signal which is a  discretization  of the radial
Ornstein-Uhlenbeck process. And, for the noise, we consider positive random variables with a
specific distribution and build a class ${\bar {\cal F}}$ such that
conditions (C1)-(C2-f)-(C3) hold.

\subsection{The signal}
 We assume that 
\begin{equation} \label{sign}
x_n= X_{n\Delta}
\end{equation}
is a discretization of a continuous time diffusion $(X_t)$ equal to a
radial Ornstein-Uhlenbeck process. We recall its definition and properties.
\subsubsection{The one-dimensional radial Ornstein-Uhlenbeck process.}
Consider the one-dimensional Ornstein-Uhlenbeck process given by:
\begin{equation}
\xi_t=  \xi_0 + \int_{0}^{t}\theta \xi_s \;ds + \sigma W_t,
\end{equation}
where $(W_t)$ is a standard Brownian motion and the initial variable
$\xi_0$ is independent of $(W_t)$. Then, 
\begin{equation} \label{ou}
\xi_t= \xi_0 \; e^{\theta t}+ \sigma \int_0^{t} e^{\theta (t-s)}\;dW_s.
\end{equation}
Let $X_t= |\xi_t|$. Then, a simple computation shows that the conditional distribution of
$X_t$ given $\xi_0=\xi$ only depends on $x=|\xi|$ so that $(X_t)$ is
a Markov process. We may call it the one-dimensional
radial Ornstein-Uhlenbeck process. Let us now give the conditional density of
$x_n=X_{n\Delta}$ given $x_{n-1}=X_{(n-1)\Delta}$, {\em i.e.} the
transition density of $(x_n)$ (see (\ref{sign})). For $\Delta >0$, we set 
\begin{equation} \label{abeta}
U_n=\xi_{n\Delta} \quad \mbox{and} \quad  a= e^{\theta \Delta}, \quad\beta^{2}= \sigma^{2} \frac{e^{2\theta\Delta}-1}{2\theta}. 
\end{equation} 
Then, as can be easily deduced from (\ref{ou}), $(U_n)$ is a standard AR($1$)-process  satisfying
\begin{equation} \label{ar1}
U_n=aU_{n-1} + \beta \eta_n, \quad n \ge 1, \quad U_0= \xi_0
\end{equation}
where  $(\eta_n, n \ge 1)$ is a sequence of
i.i.d. random variables having distribution ${\cal
  N}(0,1)$. 
Now, $x_n=|U_n|= X_{n \Delta}$ is a Markov chain having transition density (for  positive $x$)
\begin{equation} \label{tr}
p^{(1)}(x,x')= (p(x,x')+p(x,-x')) 1_{(x'>0)},
\end{equation} 
where $p(u,u')$ is the transition density of $(U_n)$, {\em i.e.}
\begin{equation} \label{trar}
p(u,u')= \frac{1}{\beta \sqrt{2\pi}} \exp{-\frac{(u'-au)^{2}}{2\beta^{2}}}.
\end{equation} 
A simple computation shows that (with $x>0$)
\begin{equation} \label{tr1}
p^{(1)}(x,x')= \;1_{(x'>0)} \frac{2}{\beta \sqrt{2\pi}}
\exp{(-\frac{x'^{2}}{2\beta^{2}})}\;\exp{(-\frac{a^{2}x^{2}}{2\beta^{2}})}
\left(\cosh{(\frac{axx'}{\beta^{2}})}\right).
\end{equation} 
Now, using the series expansion of $\cosh$, we obtain a representation
of this transition density as the following mixture of distributions:
\begin{equation} \label{tr2}
p^{(1)}(x,x')= \;1_{(x'>0)} \frac{2}{\beta \sqrt{2\pi}}
\exp{(-\frac{x'^{2}}{2\beta^{2}})}\;\sum_{i \ge 0} \alpha_i(x)
\frac{x'^{2i}}{C_{2i} \beta^{2i}},
\end{equation}
where 
$C_{2i}=\frac{(2i)!}{2^{i}i!}$ is the $2i$-th moment of a standard Gaussian variable (see
(\ref{momg})), and for $i \ge 0$,
\begin{equation} \label{tr2b}
\alpha_i(x)= \exp{(-\frac{a^{2}x^{2}}{2\beta^{2}})}
\left(\frac{a^{2}x^{2}}{2\beta^{2}}\right)^{i} \frac{1}{i!}.
\end{equation}
For $\theta<0$ ($0<a<1$), the process $(x_n)$ has a stationary
distribution given by 
\begin{equation} \label{sta}
\pi^{(1)}(dx)=  \;1_{(x>0)} \frac{2}{\rho \sqrt{2\pi}}
\exp{(-\frac{x^{2}}{2\rho^{2}})} \;dx
\end{equation}
with (see (\ref{abeta}))
\begin{equation} \label{rho}
\rho^{2}= \frac{\beta^{2}}{1-a^{2}}= \frac{\sigma^{2}}{2 |\theta|}.
\end{equation} 
These formulae will be the useful  tool for the construction of the class ${\bar
  {\cal F}}$ below.

\subsubsection{ The $\delta$-dimensional radial Ornstein-Uhlenbeck process.}
For $\delta>1$, consider the stochastic differential equation
\begin{equation} \label{radialou}
dX_t= (\theta X_t +  \frac{\sigma^{2}(\delta-1)}{2 X_t}) dt+
\sigma d\beta_t,  \quad X_0= \eta.
\end{equation}
where  $(\beta_t)$ is a
standard Brownian motion and $\eta$ is a random variable independent
of $(\beta_t)$. The values $\theta, \sigma, \delta$ are constant
parameters. This process is called a radial Ornstein-Uhlenbeck
process. This is due to the fact that, when $\delta$ is an integer
greater than $1$, then $X_t$ is the Euclidian norm of a $\delta$-dimensional vector
$(\xi_t^{1}, \ldots, \xi_t^{\delta})$ 
whose components are i.i.d. Ornstein-Uhlenbeck processes satisfying:
$$
d\xi_t^{j}= \theta \xi_t^{j} dt+
\sigma dW_t^{j}.
$$
Moreover, when $(X_t)$ is given by (\ref{radialou}),  the process
$R_t=X_t^{2}$ is the classical Cox-Ingersoll-Ross  diffusion model given by:
\begin{equation} \label{cir}
dR_t= (2 \theta R_t+ \delta \sigma^{2}) dt + 2 \sigma R_t^{1/2} d\beta_t.
\end{equation}
The processes $(X_t)$   and $(R_t)$ have explicit transition
probabilities with densities with respect to the Lebesgue measure on
$(0, +\infty)$. There are closed-form formulae for the transition
densities when $\delta$ is an odd integer. Otherwise, they
depend on Bessel functions and have  explicit developments
as sums of series.
As above, we set $x_n=X_{n \Delta}$ and give the expression of the
transition operator of this Markov chain using the notations
(\ref{abeta}). For $\alpha \ge 0$, let us set
\begin{equation} \label{calf}
C_{\alpha}= \E|X|^{\alpha},
\end{equation}
for $X$ a standard Gaussian random variable.
\begin{prop} \label{traou} 
\begin{enumerate}
\item Assume $\delta>1$. Then the transition density of the
  $\delta$-dimensional radial Ornstein-Uhlenbeck process (see
  (\ref{radialou})) is equal to (with $x>0$):
\begin{equation} \label{trd}
p^{(\delta)}(x,x')=  \;1_{(x'>0)} \frac{2}{\beta \sqrt{2\pi}}
\exp{(-\frac{x'^{2}}{2\beta^{2}})}\;\sum_{k \ge 0} \alpha_k(x)
\frac{x'^{\delta-1+2k}}{C_{\delta-1+2k} \; \beta^{\delta-1+2k}},
\end{equation}
where the mixture coefficients  are given by (\ref{tr2b}) and the couple $(a,\beta^{2})$ is linked with the original
parameters $(\theta, \sigma^{2})$ through relations (\ref{abeta}).
\item Assume that $\delta= 2n+1$ with $n \ge 1$  an 
  integer. Define the operator ${\cal T}$, acting on functions $f \in
  C^{1}((0,+\infty), \R)$,  by ${\cal T}(f)(x)=
  \frac{f'(x)}{x}$. Then, the transition density of
  $(x_n=X_{n\Delta})$, where $(X_t)$ is the $2n+1$-dimensional radial
  Ornstein-Uhlenbeck process is equal to (with $x>0$)
\begin{equation} \label{trdn}
p^{(2n+1)}(x,x')=  \;1_{(x'>0)} \frac{2}{\beta \sqrt{2\pi}}
\exp{(-\frac{x'^{2}}{2\beta^{2}})}\;\exp{(-\frac{a^{2}x^{2}}{2\beta^{2}})} \;\frac{x'^{2n}}{\beta^{2n}}
\;[{\cal T}^{n}(\cosh)(z)]_{z=\frac{axx'}{\beta^{2}}}
\end{equation}
where ${\cal T}^{n}= {\cal T} \circ {\cal T} \ldots \circ {\cal T}$ is the $n$-th iterate of
${\cal T}$.

Moreover, the above formulae also hold for $\delta=1$ ($n=0$) (see
(\ref{tr1})-(\ref{tr2})) with the convention that $T^{0}(f)= f$.
\end{enumerate}
\end{prop}

When $\theta<0$ ($0<a<1$), the process $(X_t)$ and the Markov chain $(x_n)$ have a stationary
distribution equal to (see (\ref{abeta}))-(\ref{rho}))
\begin{equation} \label{stan}
\pi^{(\delta)}(dx)=  \;1_{(x>0)} \frac{2}{\rho \sqrt{2\pi}}
\exp{(-\frac{x^{2}}{2\rho^{2}})} \frac{x^{\delta-1}}{\rho^{\delta-1}\;C_{\delta-1}} \;dx.
\end{equation}

Details are given in the Appendix.

\subsection{Distribution of the noise}
Assume that, for all $n$, $w_n$ has the distribution of $\Gamma^{-1/2}$
where $\Gamma$ has an exponential  distribution with parameter $\lambda>0$. 
Then, for all  positive $x$, the distribution of $Y=x\,w_1$  is given by:
\begin{equation} \label{dy/x}
F_{x}(dy)= f_x(y)\;dy, \quad \mbox{with} \quad f_x(y)=\frac{2 \lambda x^{2 }}{y^{3}} 
\exp(-\frac{\lambda x^2}{y^2}) \;1_{(0,\infty)}(y).
\end{equation} 
It is worth noting that  the  distribution  of $w_1$ ($F_{1}(dy)$) satisfies:
\begin{equation}
\E(|\log{w_1}|) <\infty \quad \mbox{and} \quad \E(w_1^{r}) <\infty
\quad \mbox{if and only if} \quad r< 2.
\end{equation}
Instead of an exponential distribution, we could take a Gamma
distribution with integer index.

\subsection{The class of distributions.}

First, for $i \ge 0$, we set, for $\delta \ge 1$, 
\begin{equation} \label{gi}
g_i^{(\delta)}(x)= \;1_{(x>0)}\; \frac{2}{(2\pi)^{1/2}} \frac{x^{\delta-1+2i}}{C_{\delta-1+2i}}
 \exp(-\frac{x^2}{2 }),
\end{equation}
where  $C_{\alpha}$ is defined in (\ref{calf}). Thus,
$g_i^{(\delta)}$ is a  probability density on $(0,+\infty)$.  Then, for $\sigma>0$, we set
\begin{equation} \label{nuisig}
\nu^{i,(\delta)}_{\sigma}(dx)= \frac{1}{\sigma} g_i^{(\delta)}(\frac{x}{\sigma}) dx= \;1_{(x>0)}\;\frac{2}{\sigma(2\pi)^{1/2}}  \frac{x^{\delta-1+2i}}{C_{\delta-1+2i}\,\sigma^{\delta-1+2i}}
 \exp(-\frac{x^2}{2 \sigma^{2} })\;dx.
\end{equation}
For each $i$,
the distribution $\nu^{i,(\delta)}_{\sigma}$ is equal to the  distribution of
$\sqrt{G_i^{(\delta)}}$ where  $G_i^{(\delta)}$ is  Gamma with parameters
$(i+\delta/2,1/2\sigma^2)$. This  Gamma distribution is identical to a
$\sigma^{2} \chi^{2}(\delta+2i)$ (with non integer parameter). Hence, as $i$ increases, the distributions
$(\nu^{i,(\delta)}_{\sigma}, i \ge 0)$ are stochastically increasing.
Let us now define
\begin{equation} \label{classfn}
{\cal F}^{i,(\delta)}= \{ \nu^{i,(\delta)}_{\sigma}; \sigma >0 \}, \quad {\cal F}^{(\delta)}= \cup_{i \ge 0} {\cal F}^{i,(\delta)}.
\end{equation}
And, 
\begin{equation} \label{classfnbar}
 {\bar {\cal F}^{(\delta)}}= \{\nu= \nu_{\sigma, \alpha}= \sum_{i \ge 0} \alpha_i\;
\nu^{i,(\delta)}_{\sigma}, \alpha=(\alpha_i) \in S,   \nu^{i,(\delta)}_{\sigma} \in {\cal
  F}^{i,(\delta)}, i \ge 0 \}.
\end{equation}
The  law of
$|X|$ for $X$ a Gaussian variable  with mean $m$ and variance
$\sigma^2$ has the density
\begin{equation} \label{abs}
g(x)= 1_{(x>0)} \frac{2}{\sigma \sqrt{2\pi}}
\;\exp{(-\frac{x^{2}}{2\sigma^{2}})}\exp{(-\frac{m^{2}}{2 \sigma^{2}})}
\left(\cosh{(\frac{m\,x}{\sigma^{2}})}\right).
\end{equation} 
A Taylor series development of the $\cosh$ yields a distribution of
the class  $ {\bar {\cal F}^{(1)}}$ with mixture parameter 
\begin{equation} \label{mixture}
\alpha_i(m,\sigma)=  \exp{(-\frac{m^{2}}{2 \sigma^{2}})} \frac{1}{i!}
\left(\frac{m^{2}}{2\sigma^{2}}\right)^{i}, \quad i\ge 0 .
\end{equation}
All  distributions in $ {\bar {\cal F}^{(\delta)}}$  have a density with
respect to a $\nu^{0,(\delta)}_{\sigma}(dx)$ for some positive $\sigma$ (see
(\ref{nuisig})). And this density is expressed as an entire series
of even powers. The transition density  (\ref{tr1}) and the stationary
density (\ref{sta}) belong to  ${\bar {\cal F}^{(1)}}$. The transition density (\ref{trd})
and the stationary density (\ref{stan}) belong to
$ {\bar {\cal F}^{(\delta)}}$. Thus,  our condition (C3)
holds for the model defined by (\ref{sign}) and (\ref{radialou}).

\subsection{Up-dating, marginal and prediction operators.}
In this section, we show that the filtering and prediction algorithms evolve in the
class $ {\bar {\cal F}^{(\delta)}}$ when the signal is a discrete regular sampling
of the $\delta$-dimensional radial Ornstein-Uhlenbeck process. For
this, it is enough to check conditions (C1)-(C2)-(C3).  

We have already noted that condition (C3) holds. Actually, the form of the transition
density (\ref{trd}) and of the stationary distribution (\ref{stan})
(when it exists) indicates how to define the class $ {\bar {\cal F}^{(\delta)}}$.

\begin{prop} \label{updateou}
Let $y>0$, and consider the up-dating operator $\varphi_y$ (see (\ref{phi}))
corresponding to $f_x(y)$ given in (\ref{dy/x}). For $i \ge 0$ and
  $\sigma>0$, with $\nu^{i,(\delta)}_{\sigma}$ defined in (\ref{nuisig}), 
$$
\varphi(\nu^{i,(\delta)}_{\sigma})= \nu_{T_{y}(\sigma)}^{i+1,(\delta)}, \quad \mbox{with}
\quad T_{y}(\sigma)= \frac{\sigma \; y}{(y^{2}+ 2 \lambda \sigma^{2})^{1/2}}.
$$
Thus, $t_{y}(i)=t(i)=i+1$  and (C1) holds for
$(\varphi_{y}, {\cal F}^{(\delta)})$.
\end{prop}

\proof We have 
$$
f_x(y)g_i^{(\delta)}(x/\sigma) \propto x^{\delta-1+2(i+1)} \exp{-((\frac{2\lambda}{y^{2}} +
  \frac{1}{\sigma^{2}})\frac{x^2}{2})}.
$$

Hence, we may define
$$
\frac{1}{(T_y(\sigma))^{2}}= \frac{2\lambda}{y^{2}} + \frac{1}{\sigma^{2}}.
$$
 Moreover, we have $t_{y}(i)=t(i)=i+1$: if $\nu
\in {\cal F}^{i,(\delta)}$, then, $\varphi_{y}(\nu) \in {\cal F}^{i+1,(\delta)}$. So, we
get the result.
\endproof

We also give the marginal distribution.
\begin{prop} \label{margiou}
For $i \ge 0$ and
  $\sigma>0$, with $\nu^{i,(\delta)}_{\sigma}$ defined in (\ref{nuisig}), the
  marginal density  (see (\ref{marg})) is equal to
\begin{equation}
p_{\nu^{i,(\delta)}_{\sigma}}(y)=   \; 1_{y>0}\frac{2}{\lambda^{1/2}\sigma}p_i^{(\delta)}(\frac{y}{\lambda^{1/2}\sigma}) 
\end{equation}
with
\begin{equation}
p_i^{(\delta)}(y)= \frac{(\delta+2i)y^{\delta-1+2i}}{(y^2+ 2)^{i+1+\delta/2}}.
\end{equation}
\end{prop}
\proof
Using $C_{\delta-1+2(i+1)}=(\delta+2i)C_{\delta-1+2i}$, 
we get
$$
\int_{\R^{+}}\;\frac{1}{\sigma} g_i^{(\delta)}(\frac{x}{\sigma})\;f_x(y) dx= 
\frac{\lambda \sigma^2 (\delta + 2i)y^{\delta-1+2i}}{(y^2+ 2 \lambda \sigma^2)^{i+1+\delta/2}}=
\lambda^{-1/2}\sigma^{-1} p_i^{(\delta)}(\frac{y}{\lambda^{1/2}\sigma}).
$$

\endproof

\noindent
{\bf Remark.}  Proposition \ref{margiou} allows to obtain the  density of
  $(y_0, y_1, \ldots,y_n)$, {\em i.e.} the exact likelihood based on
  this observation. Indeed, this joint density is obtained as the
  product of the conditional densities of $y_i$ given $y_{i-1},
  \ldots, y_1$. These are computed as marginal densities (see (\ref{margn})).

\begin{prop} \label{prediou}
Consider the transition operator  $P^{(\delta)}$ with transition density
(\ref{trd}). Let $\nu=\nu^{i,(\delta)}_{\sigma}$ be given by (\ref{nuisig})
with $i \ge 0$, {\em i.e.} $\nu^{i,(\delta)}_{\sigma} \in  {\cal F}^{i,(\delta)}$.  Set $\nu P^{(\delta)}= \psi^{(\delta)}(\nu)$. Then, 
\begin{equation}
\psi^{(\delta)}(\nu^{i,(\delta)}_{\sigma})= \sum_{k=0}^{i} \alpha_k^{(i,\sigma)} \nu^{k,(\delta)}_{\tau(\sigma)}
\end{equation}
with 
\begin{equation}
\tau^{2}(\sigma)= \beta^{2}+ a^{2} \sigma^{2},
\end{equation}
and for $k=0,1, \ldots,i$,
\begin{equation}
\alpha_k^{(i,\sigma)}=\binom{i}{k} \left(1 -
  \frac{\beta^2}{\tau^{2}(\sigma)}\right)^{k} \left(\frac{\beta^2}{\tau^{2}(\sigma)}\right)^{i-k} 
\end{equation}
Thus, $L(i)=i$  and $\psi^{(\delta)}(\nu^{i,(\delta)}_{\sigma})$ belongs to ${\bar
  {\cal F}}^{(\delta)}$. Condition (C2-f) holds for $(P^{(\delta)}, {\bar
  {\cal F}}^{(\delta)})$. 
\end{prop}
\proof
We have to compute 
\begin{equation} \label{nuPn}
A=\int_{0}^{\infty}  (1/\sigma)g_{i}^{(\delta)}(x/\sigma) p^{(\delta)}(x,x')\;dx,
\end{equation}
with $g_i^{(\delta)}$ given in (\ref{gi}) and $p^{(\delta)}(x,x')$ given in
(\ref{trd}). 
Let us define $s^2$ by
\begin{equation} \label{s}
\frac{1}{s^2}= \frac{a^2}{ \beta^2}+ \frac{1}{\sigma^2}= \frac{\tau^{2}(\sigma)}{\beta^{2}\sigma^{2}}.
\end{equation}
Hence, for all $k \ge 0$, 
\begin{equation}
\int_{0}^{\infty} 2 \;x^{2(i+k)+\delta-1}
\exp{(-\frac{a^{2}x^{2}}{2\beta^{2}})}\exp{(-\frac{x^{2}}{2\sigma^{2}})}
 \frac{dx}{ (2\pi)^{1/2}}= C_{2(i+k)+\delta-1} \; s^{2(i+k)+\delta}.
\end{equation}
Now, using $s/\beta \sigma= 1/\tau(\sigma)$,  $A$ is given as the following expression
\begin{equation} \label{A}
A= \frac{2\; 1_{(x'>0)}}{\tau(\sigma) (2\pi)^{1/2}}
\exp{(-\frac{x'^{2}}{2\beta^{2}})} \; \frac{x'^{\delta-1}}{\tau^{\delta-1}(\sigma)} \;\Sigma,
\end{equation}
with
\begin{equation} \label{Sig}
\Sigma= \sum_{k=0}^{\infty}  \frac{x'^{2k}}{\beta^{2k}}
\left(\frac{a^{2}}{2
    \;\beta^{2}}\right)^{k}\frac{C_{2(i+k)+\delta-1}}{k!C_{2k+\delta-1}C_{2i+\delta-1}} \frac{s^{2(i+k)}}{\sigma^{2i}}.
\end{equation}
Using  (\ref{s}) and some computations, we get, for all $k \ge 0$,  
\begin{equation}
\left(\frac{a^{2}}{\beta^{2}}\right)^{k}\frac{s^{2(i+k)}}{\sigma^{2i}}=\left(\frac{\beta^{2}}{\tau^{2}(\sigma)}\right)^{i}\left(1-\frac{\beta^{2}}{\tau^{2}(\sigma)}\right)^{k}.
\end{equation}
Now, we set 
\begin{equation} \label{c}
\frac{1}{c^2}= \frac{1}{\beta^2} \left(1- \frac{\beta^{2}}{\tau^{2}(\sigma)}\right).
\end{equation}
This yields
\begin{equation} \label{Sig2}
\Sigma= \sum_{k=0}^{\infty} \frac{1}{k!}\left( \frac{x'^{2}}{2 \;c^{2}}\right)^{k}
\left(\frac{\beta^{2}}{\tau^{2}(\sigma)}\right)^{i}\frac{C_{2(i+k)+\delta-1}}{C_{2k+\delta-1}C_{2i+\delta-1}} .
\end{equation}
Now, we use the following lemma whose proof is given in the Appendix.
\begin{lemma} \label{lem}
For all integer $k \ge 0$, and all $i \ge n$
\begin{multline}
\frac{C_{2(i+k)+\delta-1}}{C_{2k+\delta-1}C_{2i+\delta-1}}= \sum_{j=0}^{i}
k(k-1)\ldots(k-j+1) a_j\\
= a_0+ k a_1 + k(k-1) a_2 + \ldots + k(k-1) \ldots (k-i+1) a_{i}.
\end{multline}
with, for $j=0,1, \ldots, i$, 
$$
a_j= \binom{i}{j} \frac{2^{j}}{C_{2j+\delta-1}},
$$
and the coefficient of $a_0$ is equal to $1$.
\end{lemma}
Now, we  transform expression (\ref{Sig2}) into
\begin{equation}
\Sigma= \sum_{j=0}^{i} \Sigma_j,
\end{equation}
where
\begin{equation}
\Sigma_j=a_j
\left(\frac{\beta^{2}}{\tau^{2}(\sigma)}\right)^{i}\sum_{k=0}^{\infty}
\frac{1}{k!}\left( \frac{x'^{2}}{2 \;c^{2}}\right)^{k} \;k(k-1)\ldots(k-j+1).
\end{equation}
But, $k(k-1)\ldots(k-j+1)=0$ for $k=0,1, \ldots, j-1$. So, we get
\begin{multline}
\Sigma_j=a_j \left(\frac{\beta^{2}}{\tau^{2}(\sigma)}\right)^{i}\sum_{k \ge j}
\left( \frac{x'^{2}}{2 \;c^{2}}\right)^{k}\frac{1}{(k-j)!}\\= 
\frac{1}{C_{2j+\delta-1}} \left(\frac{\beta^{2}}{\tau^{2}(\sigma)}\right)^{i} \binom{i}{j} \left( \frac{x'^{2}}{ \;c^{2}}\right)^{j}\exp{(\frac{x'^{2}}{2\;c^{2}})}.
\end{multline}
Now, we compute $A$ from (\ref{A}). Using (\ref{c}), we add the
exponents of the
exponential terms and after some elementary computations, we obtain
\begin{equation}
A= \frac{2\; 1_{(x'>0)}}{\tau(\sigma) (2\pi)^{1/2}}
\exp{(-\frac{x'^{2}}{2\tau^{2}(\sigma)})} \; \Sigma'
\end{equation}
with 
\begin{equation}
\Sigma'= \sum_{j=0}^{i} \Sigma'_j,
\end{equation}
and 
\begin{equation}
\Sigma'_j= \binom{i}{j} \frac{1}{C_{2j+\delta-1}}\left( \frac{x'^{2j+\delta-1}}{\tau^{2j+\delta-1}(\sigma)}\right)\left(\frac{\beta^{2}}{\tau^{2}(\sigma)}\right)^{i-j}\left(1-\frac{\beta^{2}}{\tau^{2}(\sigma)}\right)^{j}.
\end{equation}
So the proof is complete.
\endproof

\noindent
{\bf Remarks.}
\begin{enumerate}
\item Here, we have two representations of
  $\psi^{(\delta)}(\nu^{i,(\delta)}_{\sigma})$. One has scale parameter
  $\beta$ and a mixture parameter with infinite length. The second has
  scale parameter $\tau(\sigma)$ and a finite length mixture
  parameter. The latter appears as a minimal representation of this
  distribution in a sense that we try to clarify (work in progress).
\item By Propositions  \ref{updateou} and \ref{prediou}, we see that
  both the exact filter and the prediction filter evolve in the
  extended class ${\bar {\cal F}}^{(\delta)}$. If the initial distribution
  of the signal is in the subclass ${\cal F}^{(\delta)}_{f}$  ({\em e.g.} if the
  signal is in stationary regime),   it has a mixture coefficient
  with finite length (see (\ref{classfnbar})). Then, the number of
  mixture components grows of a unit at each iteration but remains finite. However, the
  numerical simulations that we have done in Genon-Catalot and
  Kessler (2004) show that there are only two or three significantly
  non nul mixture coefficients. Some stability results are also
  obtained that may be extended to the model investigated here. 
\end{enumerate}

\section{Stochastic volatility type models.} \label{scalecir}
We first draw some immediate consequences of the previous
section. Then, we  introduce some new type of stochastic
volatility models.
\subsection{Scale perturbation of a Cox-Ingersoll-Ross diffusion process.}
Consider now the model obtained by taking squares of the previous
one. Set 
\begin{equation}
z_n= y_n^{2}=r_n v_n ,
\end{equation}
with $r_n=x_n^{2}$,  $v_n=w_n^{2}$ and $(x_n,w_n)$ as in the previous
section. Then $r_n=  R_{n\Delta}$ is a  discrete sampling of the Cox-Ingersoll-Ross diffusion model
(\ref{cir}). It is a  Markov chain with transition (see (\ref{trd}))
\begin{multline}
q^{(\delta)}(r,r')=  (1/2) r'^{-\frac{1}{2}}\;p^{(\delta)}(r^{1/2},r'^{1/2} )\\
=1_{(r'>0)} \frac{1}{\beta \sqrt{2\pi}}
\exp{(-\frac{r'}{2\beta^{2}})}\;\sum_{k \ge 0} \alpha_{k}(r^{\frac{1}{2}})
\frac{(r')^{k-1+\frac{\delta}{2}}}{C_{2k+\delta-1} \beta^{2k+\delta-1}}.
 \end{multline}
(see (\ref{tr2b}) for $(\alpha_k(x))$.
This is now a mixture of Gamma distributions with parameters
$(k+\frac{\delta}{2}, \frac{1}{2 \beta^{2}})$. 

The distribution of the noise $(v_n)$ is now inverse exponential. And
the class of distributions is composed with mixtures of Gamma
distributions with parameters $(i+\frac{\delta}{2}, \frac{1}{2
  \sigma^{2}})$, for  $i \ge 0$. The filtering and prediction
algorithm can be explicitely expressed with the same formulae after
some simple changes  $z=y^{2}$ for the observations and the change of
variables 
$x=r^{\frac{1}{2}}$ for the distributions.

Note that another computable filter is obtained by setting
\begin{equation} \label{inverse}
z'_n= \frac{1}{z_n}= \frac{1}{r_n} v'_n, \quad v'_n= \frac{1}{v_n}.
\end{equation}
\subsection{Stochastic volatility type models.}
The above considerations lead to some new type of stochastic
volatility models. Indeed, stochastic volatility models usually
postulate that the observed price process  $(S_n)$ of an asset is such that
\begin{equation} \label{svm}
Z_n= \log{\frac{S_{n+1}}{S_n}}= \sqrt{V_n} \varepsilon_n,
\end{equation}
where $(V_n)$ is a positive Markov chain (the unobserved volatility),
$(\varepsilon_n)$ is a sequence of i.i.d. standard Gaussian variables,
the two sequences being independent. We do not know explicit filters for such stochastic volatility models when the signal is a
discrete sampling of a diffusion process. 

Now, taking squares in (\ref{svm}), we get that $Z_n^{2}=V_n
\varepsilon_n^{2}$, where $\varepsilon_n^{2}$ is distributed as a
$\chi^{2}(1)=G(1/2, 1.2)$. Our previous study suggests to replace the
$G(1/2, 1/2)$ distribution by a $G(1,\lambda)$ (possibly a
$G(k,\lambda)$ with $k$ integer). More precisely, the following
stochastic volatility type models will provide explicit filters
through a symetrization device.  Consider 
\begin{equation} \label{svmt}
Z'_n=  \sqrt{r_n} \varepsilon'_n, \quad \mbox{or} \quad Z''_n= \frac{1}{\sqrt{r_n}} \varepsilon''_n,
\end{equation}
with $r_n=R_{n\Delta}$ a discrete sampling of a Cox-Ingersoll-Ross
diffusion. For the noises, consider a symetric Bernoulli variable
$\varepsilon \pm 1$ with probability $1/2$, independent of a random
variable $\Gamma$ having distribution $G(1,\lambda)$ (exponential
distribution). Then, assume that $\varepsilon'_n$ is distributed as
$\frac{\varepsilon}{\sqrt{\Gamma}}$ and that $\varepsilon''_n$ is distributed
as $\varepsilon\sqrt{\Gamma}$. Then, the two models of filtering given in
(\ref{svmt}) can be solved explicitely.

\section{ A discretized Cox-Ingersoll-Ross diffusion and  conditionally Poisson observations.} \label{poisson}
Models which have no  representation  as $y_n=H(x_n, w_n)$, for some simple
 function $H$, are also of
interest. These models are completely specified by the conditional
distribution (\ref{condit}) and the transition operator of the hidden
Markov chain. We investigate below such an  example. Suppose that the couple (signal, observation) is defined as follows.
The signal is the process $(r_n)$ obtained as above from a discrete
sampling of the square-root model (\ref{cir}). Now, the observation
$y_n$ is such that, given $r_n=r$, $y_n$ has a Poisson distribution
with parameter $\lambda \;r$, {\em i.e.}
\begin{equation}
\P(y_n=y|r_n=r) = f_r(y)= \exp{(- \lambda r)} \frac{(\lambda\;r)^{y}}{y!}, y \in \N.
\end{equation}
We can check our conditions with the class of distributions fitted with
the signal $(r_n)$, i.e. the class of mixtures of Gamma distributions $ G(i+\frac{\delta}{2}, \frac{1}{2 \sigma^{2}})$
with parameters  $(i+\frac{\delta}{2}, \frac{1}{2 \sigma^{2}})$, for  $i \ge 0$. 
Only (C1) needs to be checked.  Let us set for this Gamma density
\begin{equation}
\gamma_{\sigma}^{i,(\delta)}(r)=1_{(r>0)} \frac{1}{\beta \sqrt{2\pi}}
\exp{(-\frac{r}{2\sigma^{2}})}
\frac{r^{i-1+\frac{\delta}{2}}}{C_{2i+\delta-1} \sigma^{2i+\delta-1}}.
 \end{equation}
Now, 
\begin{equation}
f_r(y) \; \gamma_{\sigma}^{i}(r) \propto 
\exp{[-(\lambda+\frac{1}{2\sigma^{2}})r]}\; r^{(y+ i-1+\frac{\delta}{2})}.
\end{equation}
This is again a Gamma distribution of the same type. And we get
\begin{equation}
t_{y}(i)= y+i,\quad  T_{y}(\sigma)= (2\lambda + \frac{1}{\sigma^{2}})^{-\frac{1}{2}}.
\end{equation}
The filtering and prediction algorithm will
evolve in the family of mixtures of Gamma distributions with tail index
$i+(\delta/2)$, $i \ge 0$.
The marginal distributions can be explicitely computed.

\section{Concluding  remarks.} \label{concl}
This work has to be completed by numerical simulations. In Genon-Catalot and Kessler
(2004), the model corresponding to a one-dimensional radial
Ornstein-Uhlenbeck process is studied and implemented. The numerical
results show that  the number of significantly non
nul mixture coefficients is less than $2$ or $3$. Theoretical properties
linked with the stability of the filters are established in this paper
which may be
extended to the models of Section \ref{scalepert}. 

The above results may be  extended to the case of a
non time-homogeneous signal. For instance, it is possible to  consider a non regular discrete sampling
of the underlying diffusion model ({\em i.e.} to consider
$x_n=X_{t_{n}}$ with $0<t_1<\ldots <t_n<\ldots$). It is also possible to
consider non time-homogeneous conditional distributions of the
observation given the signal. Note also that the signal may or may not
be ergodic.

\noindent
{\bf Acknowledgments.} The authors wish to thank Wolfgang Runggalddier
and Pavel Chigansky  for helpful discussions and references.

\section{Appendix}
\subsection{The extended Kalman filter}
\noindent
{\em The class ${\cal F}^{i}$}.  We compute the Laplace transform of a
distribution  $\nu^{i}_{(\mu,m,\sigma^{2})}$ (see (\ref{gnui})). For
$X$ a random variable having the previous distribution, 
elementary computations using (\ref{cst1})-(\ref{cst2}) yield, for
$\lambda \in \mathbb{C}$,
\begin{equation} \label{laplace}
 E(\exp{\lambda X})= \;\frac{C_{2i}(m+\mu+\lambda
  \sigma^{2};\sigma^{2})}{C_{2i}(m+\mu;\sigma^{2})}\; \exp{(\lambda m+ \frac{\lambda^{2}}{2
    \sigma^{2}})} .
\end{equation}
All parameters $(i,m,\mu,\sigma^{2})$  are identifiable.
>From this formula, we can prove that, for all $m,\mu$ and all $i$, as $\sigma$ tends to $0$,
$\nu^{i}_{(\mu,m,\sigma^{2})}$ weakly converges to the Dirac measure
$\delta_{m}$. Moreover, for any mixture coefficient $\alpha$, $\sum_{i
  \ge 0} \alpha_i \;\nu^{i}_{(\mu,m,\sigma^{2})}$ weakly converges also
to $\delta_m$.

\noindent
{\em Proof of Proposition \ref{upd}}. Let us consider a random
variable $X$ with distribution $\nu=\nu^{i}_{(\mu,m,\sigma^{2})}$ and
let $Y=hX + \varepsilon$ with $X$ and $\varepsilon$ independent, and
$\varepsilon$ having distribution ${\cal N}(0, \gamma^{2})$. Then,
$\varphi_{y}(\nu)$ is exactly the conditional distribution of $X$
given $Y=y$.  Its density is proportional to:
\begin{equation} \label{de}
x \rightarrow (x+\mu)^{i} \; \exp{[-(\frac{(y-hx)^{2}}{2 \gamma^{2}}+\frac{(x-m)^{2}}{2 \sigma^{2}})]}
\end{equation}
We compute the exponent of the exponential above and obtain:
\begin{equation}
\frac{(y-hx)^{2}}{2 \gamma^{2}}+\frac{(x-m)^{2}}{2 \sigma^{2}}=
\frac{(x-{\hat m}(y))^{2}}{2 {\hat \sigma}^{2}},
\end{equation}
with:
\begin{equation} \label{hatbar}
{\hat m}(y)= \frac{m \gamma^{2}+h y \sigma^{2}}{{\bar \sigma}^{2}},
\quad {\bar \sigma}^{2}=\gamma^{2}+ h^{2} \sigma^{2},
\end{equation}
and
\begin{equation} \label{hatsigm}
{\hat \sigma}^{2}= \frac{\sigma^{2} \gamma^{2}}{{\bar \sigma}^{2}}.
\end{equation}
This implies that $\varphi_{y}(\nu)=\nu^{i}_{(\mu,{\hat m}(y),{\hat
    \sigma}^{2})}$. So, we get the proposition. Note that this result
contains the standard case where $i=0$ and $\mu=0$.

\noindent
{\em Proof of Proposition \ref{pre}}. As noted in the text above, in this model, the transition kernel $p(x,x')$ and the
conditional  kernel $f_{x}(y)$ are of the same form. Therefore, the
computations of $\psi( \nu)$ and of the marginal density
$p_{\nu}(y)$ of $Y$ are identical up to a change of notations
($(a,\beta^{2})$ for $\psi( \nu)$, and $(h, \gamma^{2})$ for
$p_{\nu}(y)$). Because of the previous proof, it is more convenient  
here to  compute the marginal density of $Y$ when $Y=h X+
\varepsilon$ and $(X,\varepsilon)$ are as in the previous proof. We
shall use the same notations  as in the statement of Proposition
\ref{pre}, but the formulae will be given with $(h,\gamma^{2})$. We have to
integrate $f_{x}(y) \; \nu(dx)$
(with $\nu=\nu^{i}_{(\mu,m,\sigma^{2})}$) with respect to $x$. After
some elementary computations, we obtain:
\begin{equation}
p_{\nu}(y) = A_{2i}\; \exp{[-\frac{(y-{\bar m})^{2}}{2
   {\bar  \sigma}^{2}}]}, \quad {\bar m}=hm,
\end{equation}
where (see (\ref{hatbar})-(\ref{hatsigm}))
\begin{equation}
A_{2i}= \frac{C_{2i}(\mu+ {\hat m}(y);{\hat \sigma}^{2})}{C_{2i}(\mu+m;\sigma^{2})}.
\end{equation}
Let us set (see (\ref{bi}))
\begin{equation}
B_{i}= \sum_{k=0}^{i} \binom{i}{k}  \frac{(\mu+m)^{2k}}{C_{2k}\;\sigma^{2k}},
\end{equation}
\begin{equation}
{\bar \mu}=  \frac{m \gamma^{2} + \mu {\bar \sigma}^{2}}{h\sigma^{2}}.
\end{equation}
Thus, 
\begin{equation}
\mu + {\hat m}(y)= \frac{h \sigma^{2}}{\bar \sigma^{2}} \;(y+{\bar
  \mu}), \quad \frac{{\bar \mu}+{\bar m}}{{\bar
    \sigma}}=\frac{\mu+m}{\sigma} \frac{{\bar \sigma}}{h \sigma}.
\end{equation}
After some computations, we obtain
\begin{equation}
A_{2i}= \frac{1}{B_i}\; \sum_{k=0}^{i} \binom{i}{k} \left(\frac{h
  \sigma^{2}}{{\bar \sigma}^{2}}\right)^{2k}\;\frac{(y+{\bar
    \mu})^{2k}}{C_{2k}\;{\hat \sigma}^{2k}}.
\end{equation}
Now, we set 
\begin{equation}
a_{k,i}=\binom{i}{k} \left(\frac{h^{2}
  \sigma^{2}}{ \gamma^{2}}\right)^{2k} \sum_{j=0}^{k} \binom{k}{j}
\frac{({\bar \mu}+{\bar m})^{2j}}{C_{2j}\;{\bar \sigma}^{2j}},
\end{equation}
and 
\begin{equation}
{\bar \alpha}_{k}^{(i)}= \frac{\gamma^{2i}}{B_i \;{\bar \sigma}^{2i}} \;a_{k,i}.
\end{equation}
Finally, for $k=0,\ldots,i$, we obtain the following mixture coefficients:
\begin{equation}
{\bar \alpha}_{k}^{(i)}= \binom{i}{k} \frac{\gamma^{2(i-k)}}{B_i} \sum_{j=0}^{k} \binom{k}{j}
\frac{( \mu+ m)^{2j}}{C_{2j}\; \sigma^{2j}} \frac{h^{2(k-j)}
  \sigma^{2(k-j)}}{{\bar \sigma}^{2(i-j)}}.
\end{equation}
And 
\begin{equation}
p_{\nu}(y)dy = \sum_{k=0}^{i} {\bar \alpha}_{k}^{(i)} \nu^{k}_{({\bar
    \mu},{\bar m},{\bar \sigma}^{2})}(dy).
\end{equation}
So the proof is complete.

\subsection{The radial Ornstein-Uhlenbeck process.}
\subsubsection{Gaussian moments, Gamma function.}
Let us set, for $\alpha \ge 0$,  and $X$ a standard Gaussian variable,
\begin{equation} \label{momabs}
C_{\alpha}= \E(|X|^{\alpha})
\end{equation}
And recall the definition of the usual Gamma function 
\begin{equation}
\Gamma(a)= \int_{0}^{+\infty} x^{a-1} e^{-x} dx, a>0.
\end{equation}
The following relations are obtained by elementary computations.
\begin{equation} \label{relations}
C_{\alpha+1}= \alpha C_{\alpha-1}, \alpha \ge 1,\quad C_{\alpha}=
\frac{\Gamma(\frac{\alpha+1}{2})}{\sqrt{2 \pi} 2^{\frac{\alpha+1}{2}}},
\alpha \ge 0, \quad\Gamma(a)= \frac{\sqrt{2 \pi}}{2^a} C_{2a-1}, a \ge
1/2.
\end{equation}
Thus, when $\alpha=2i$, $i \in \N$, {\em i.e.} $\alpha$ is an even
integer, we obtain, $C_0=1$ and for $i \ge 1$,
\begin{equation} \label{momexpl}
C_{2i}= (2i-1) C_{2(i-1)}= (2i-1)(2i-3) \ldots 5.3.1= \frac{(2i)!}{2^{i}i!}.
\end{equation}

  \subsubsection{Transition densities.}
For an integer $\delta>1$, consider  processes $(\xi_t^{1}, \ldots, \xi_t^{\delta})$ 
 satisfying for all $j$:
$$
d\xi_t^{j}= \theta \xi_t^{j} dt+
\sigma dW_t^{j}
$$
where $(W_t^{j})$ are independent Wiener processes. Let us set $R_t=
\sum_{j=1}^{\delta} \xi_t^{j})^{2},  \quad X_t= R_t^{1/2}$. By the Ito
formula, we obtain $dR_t=   \sum_{j=1}^{\delta} 2
\xi_t^{j}\;d\xi_t^{j} +  \delta \sigma^{2} dt$. By L\'evy's characterization, the process defined by
$$
\beta_t= \int_{0}^{t} \frac{\sum_{j=1}^{\delta}\xi_t^{j}\;dW_t^{j}}{X_t}
$$
is a standard Brownian motion. And, 
\begin{equation} \label{cir1}
dR_t= (2 \theta R_t+ \delta \sigma^{2}) dt + 2 \sigma R_t^{1/2} d\beta_t.
\end{equation}
Therefore, the process $(R_t)$ is the classical Cox-Ingersoll-Ross diffusion
process. Another application of the Ito formula gives the stochastic
differential of $(X_t)$:
\begin{equation} \label{radou1}
dX_t= (\theta X_t+ \frac{(\delta-1) \sigma^{2}}{2 X_t}) dt +  \sigma  d\beta_t.
\end{equation}
Now, we do not assume any more that $\delta$ is an integer. We assume in
the stochastic differential equations (\ref{cir1}) and (\ref{radou1})
that $\delta$ is a real parameter satisfying $\delta>1$ and define the index
$\nu= (\delta/2)-1$. When $\theta=0$ and $\sigma=1$, $(X_t)$ is the
standard Bessel process with index $\nu$.  The scale and speed
densities of $(X_t)$ are obtained by the classical formulae for
one-dimensional diffusion processes. The scale density is given by
\begin{equation}
s(x)= \exp{(-\frac{2}{\sigma^{2}} \int^{x} (\theta u + \frac{(\delta-1) \sigma^{2}}{2 u} ) du
} \propto   x^{-(\delta-1)}\;\exp{(- \frac{\theta x^{2}}{\sigma^{2}})}.
\end{equation}
The speed density is $m(x)= s^{-1}(x)$. The diffusion process (\ref{radou1}) is positive recurrent on
$(0,+\infty)$ for
$\theta<0$. In this case, its stationary density is obtained by
normalizing $m$ into a probability density. Setting 
$$
\rho= \frac{\sigma}{(2|\theta|)^{1/2}},
$$
we obtain the stationary density
\begin{equation}
\pi^{(\delta)}(x)= \; 1_{x>0} \frac{2}{\rho \sqrt{2 \pi}}
\exp{(-\frac{x^{2}}{2 \rho^{2}})} \left(\frac{x}{\rho}\right)^{\delta-1} \frac{1}{C_{\delta-1}}.
\end{equation}
This is the distribution of $\Gamma^{1/2}$ with $\Gamma$ having Gamma
distribution $G(\delta/2,1/2 \rho^{2})$.

The processes $(X_t)$   and $(R_t)$ have explicit transition
probabilities with densities with respect to the Lebesgue measure on
$(0, +\infty)$. For these, we refer {\em e.g.} to Karlin and Taylor
(p.333-334). For the properties of Bessel functions that we use, we
refer {\em e.g.} to Nikiforov and Ouvarov (1983). The conditional density of
$X_{\Delta}$ given $X_0=x$ is as follows:
$$
p(\Delta,x,x')= p^{(\delta)}(x,x')= 2 \times\;1_{(x'>0)} (x')^{\delta-1} 
\exp{(\frac{\theta \;x'^{2}}{\sigma^{2}})} \times
$$
$$
 \exp{\left(-\frac{\theta}{\sigma^{2}} (\frac{e^{2\theta
          \Delta}}{(e^{2\theta \Delta}-1)}) (x^{2}+x'^{2})\right)} 
(\frac{\theta}{\sigma^{2}(e^{2\theta \Delta}-1)} )
(xx' e^{\theta \Delta})^{-\nu} \times
$$
$$ 
 I_{\nu}(xx' e^{\theta t}\frac{2\theta}{\sigma^{2}(e^{2\theta t}-1)} )
$$
where $I_{\nu}(z)$ is the Bessel function with index $\nu$.  This
function is given by the following series development
\begin{equation}
I_{\nu}(z)= \left(\frac{z}{2}\right)^{\nu} \sum_{k \ge 0}  \left(\frac{z}{2}\right)^{2k}
\frac{1}{k! \Gamma(k+\nu+1)},
\end{equation}
where $\Gamma$ is the usual Gamma function. Now, we use the
notations (\ref{abeta}) and the relations (\ref{relations}) to transform $p^{(\delta)}(x,x')$ and obtain:
\begin{equation}
p^{(\delta)}(x,x')=  \;1_{(x'>0)} \frac{2}{\beta \sqrt{2\pi}}
\exp{(-\frac{x'^{2}}{2\beta^{2}})}\;\sum_{k \ge 0} \alpha_k(x)
\frac{x'^{\delta-1+2k}}{C_{\delta-1+2k} \; \beta^{\delta-1+2k}},
\end{equation}
where the mixture coefficients  are given by (see (\ref{tr2b}))
$$\alpha_k(x)= \exp{(-\frac{a^{2} x^{2}}{2\beta^{2}})}\; \frac{1}{k!}\;\left(\frac{a^{2}x^{2}}{2\beta^{2}}\right)^{k}$$

 Now, when $\delta=2n+1$, the index is $\nu=n - \frac{1}{2}$  a half
 integer. Then,  the Bessel function is
explicit and equal to:
\begin{equation}
I_{n - \frac{1}{2}}(z)= \left(\frac{2}{\pi z}\right)^{1/2} \; z^{n} \; {\cal T}^{n}(\cosh(z))
\end{equation}
where ${\cal T}^{n}= {\cal T} \circ\ldots {\cal T}$is the n-th iterate
of the operator ${\cal T}(f)(x)= f'(x)/x$. And we obtain (\ref{trd}).

\subsubsection{Technical lemma.}
\noindent
{\em Proof of Lemma \ref{lem}}

Let us set
\begin{equation}
\varphi(k)= \frac{C_{2(i+k)+\delta-1}}{C_{2k+\delta-1}C_{2i+\delta-1}}= \frac{1}{C_{2i+\delta-1}}
(\delta-1+2k+2i-1)(\delta-1+2k+2i-3)\ldots (\delta-1+2k+1).
\end{equation}
Hence, $\varphi(x)$ is a polynomial of degree $i$ 
which admits a unique representation as a sum of the elementary
polynomials $1, x, x(x-1), \ldots, x(x-1) \ldots (x-i+1)$, say 
\begin{equation}
\varphi(x)= a'_0 + a'_1 x + a'_2 x(x-1) + \ldots + a'_{i} x(x-1)
\ldots (x-i+1).
\end{equation}
 Let us set
\begin{equation}
\psi(x)=  a_0 + a_1 x + a_2 x(x-1) + \ldots + a_{i} x(x-1) \ldots (x-i+1)
\end{equation}
with the coefficients $a_j$ given in the statement of Lemma \ref{lem}.
We will  prove that $\varphi$ and $\psi$ are identical. For this, it
 is enough to check that
\begin{equation} \label{phipsi}
\varphi(j)=\psi(j)= a_0 + j a_1 + j(j-1) a_2+ \ldots +j!  a_{j} \quad \mbox{for all} \quad j=0,1, \ldots, i.
\end{equation}
Computing  the constant and the higher degree terms, it is easy to
see that
$$
\varphi(0)=a_0= \frac{1}{C_{\delta-1}} \quad \mbox{and}  \quad
a'_{i}=a_{i}= \frac{2^{i}}{C_{\delta-1+2i}}.
$$
Now, let us fix $0<j<i$. Then,  $\varphi(j)$ and $\psi(j)$ have the
following expressions:
\begin{multline}
\varphi(j)=
\frac{C_{\delta-1+2(i+j)}}{C_{\delta-1+2j}C_{\delta-1+2i}}\\
= \frac{1}{C_{\delta-1+2j}}
(\delta-1+2j+2i-1)(\delta-1+2j+2i-3)\ldots (\delta-1+2i+1)\\
= P_{\delta-1}^{(j)}(i),
\end{multline}
\begin{multline}
\psi(j)= \frac{1}{C_{\delta-1}}+ j \frac{2i}{C_{\delta-1+2}}+ j(j-1)
\frac{2^{2}i(i-1)}{2 \;C_{\delta-1+4}}\\
+j(j-1)(j-2)\frac{2^{3}i(i-1)(i-2)}{3! \;C_{\delta-1+6}} + \ldots\\
+j!\frac{2^{j}i(i-1)\ldots (i-j+1)}{j! \;C_{\delta-1+2j}}= Q_{\delta-1}^{(j)}(i).
\end{multline} 
Hence, both quantities are polynomials of degree $j$ as functions of the variable
$i$. 
 
We now prove  that,  for all $j$, all $\delta$ and all $y$,
\begin{equation} \label{egal}
P_{\delta-1}^{(j)}(x)=Q_{\delta-1}^{(j)}(y),
\end{equation}
with 
\begin{equation} \label{egal0}
P_{\delta-1}^{(0)}(y)=Q_{\delta-1}^{(0)}(y)= \frac{1}{C_{\delta-1}},
\end{equation}
and
\begin{equation}
P_{\delta-1}^{(j)}(y)= \frac{1}{C_{\delta-1+2j}}(2y+\delta-1+2j-1)(2y+\delta-1+2j-3)\ldots (2y+\delta-1+1)
\end{equation}
\begin{equation}
Q_{\delta-1}^{(j)}(y)= \sum_{k=0}^{j} \binom{j}{k} 
\frac{2^{k}y(y-1)\ldots (y-k+1)}{C_{\delta-1+2k}}.
\end{equation}
Let us first look at $P_{\delta-1}^{(j)}(y)$. Using   $C_{\delta-1+2j}= (\delta-1+2j-1) C_{\delta-1+2(j-1)}$, we get
\begin{equation} \label{recur}
P_{\delta-1}^{(j)}(y)= P_{\delta-1}^{(j-1)}(y) + 2y P_{\delta+1}^{(j-1)}(y-1)
\end{equation}
Now, we look at $Q_{\delta-1}^{(j)}(y)$. Using the  relation
\begin{equation}
\binom{j}{k}= \binom{j-1}{k}+  \binom{j-1}{k-1},
\end{equation}
we obtain
\begin{equation}
Q_{\delta-1}^{(j)}(y)= Q_{\delta-1}^{(j-1)}(y) + 2y Q_{\delta+1}^{(j-1)}(y-1).
\end{equation}
Therefore, both families of polynomials satisfy the same relation (\ref{recur}). Since (\ref{egal0}) holds, we get
(\ref{egal}). So, the proof of the Lemma is now complete.

\end{document}